\documentclass[12pt]{amsart}
\usepackage{amssymb}
\usepackage{t1enc}
\usepackage[latin2]{inputenc}
\usepackage{verbatim}
\usepackage{amsmath,amsfonts,amssymb,amsthm}
\usepackage[mathcal]{eucal}
\usepackage{enumerate}
\usepackage[centertags]{amsmath}
\usepackage{graphics}

\newtheorem{theorem}{Theorem}
\newtheorem{lemma}{Lemma}

\newtheorem*{Go2}{Theorem A}

\newtheorem*{Go3}{Theorem B}

\begin{document}
\author{T. Tepnadze$^{1}$, L. E. Persson$^{2}$}
\title[Some inequalities for Ces\`{a}ro Means]{Some inequalities for Ces\`{a}%
ro Means of double Vilenkin-Fourier Series}
\date{}
\maketitle

\begin{abstract}
In this paper we state and prove some new inequalities related to the rate
of $L^{p}$ approximation by Ces\`{a}ro means of the quadratic partial sums
of double Vilenkin-Fourier series of functions from $L^{p}$.
\end{abstract}

$^{1}$ \email{tsitsinotefnadze@gmail.com} ;
The Artic University of Norway, Campus Narvik, P.O. Box 385, N-8505, Narvik, Norway.

$^{2}$The Artic University of Norway, Campus Narvik, P.O. Box 385, N-8505, Narvik, Norway.



\section{\protect\bigskip Introduction}

Let $N_{+}$ denote the set of positive integers, $N:=N_{+}\cup \{0\}.$ Let $%
m:=\left( m_{0},m_{1},...\right) $ denote a sequence of positive integers
not less then 2. Denote by $Z_{m_{k}}:=\{0,1,...,m_{k}-1\}$ the additive
group of integers modulo $m_{k}$. Define the group $G_{m}$ as the complete
direct product of the groups $Z_{m_{j}}$, with the product of the discrete
topologies of $Z_{mj}$'s.

The direct product of the measures%
\begin{equation*}
\mu _{k}\left( \{j\}\right) :=\frac{1}{m_{k}}\ \ \ \ \ \ \ \ \left( j\text{ }%
\in Z_{m_{k}}\right)
\end{equation*}%
is the Haar measure on $G_{m}$ with $\mu \left( G_{m}\right) =1.$ If the
sequence $m$ is bounded, then $G_{m}$\ is called a bounded Vilenkin group.
In this paper we will consider only bounded Vilenkin groups. The elements of 
$G_{m}$ can be represented by sequences $x:=\left(
x_{0},x_{1},...,x_{j},...\right) ,$ $\left( x_{j}\in Z_{m_{j}}\right) .$ The
group operation $+$ in $G_{m}$ is given by%
\begin{equation*}
x+y=\left( \left( x_{0}+y_{0\text{ }}\right) mod\text{ }m_{0},...,\left(
x_{k}+y_{k\text{ }}\right) mod\text{ }m_{k},...\right) ,
\end{equation*}%
where $x:=\left( x_{0},...,x_{k},...\right) $ and $y:=\left(
y_{0},...,y_{k},...\right) \in G_{m}.$ The inverse of $+$ will be denoted by 
$-.$

It is easy to give a base for the neighborhoods of $G_{m}:$

\begin{equation*}
I_{0}\left( x\right) :=G_{m},
\end{equation*}

\begin{equation*}
I_{n}\left( x\right) :=\{y\in G_{m}|y_{0\text{ }}=x_{0\text{ }},...,y_{n-1%
\text{ }}=x_{n-1\text{ }}\}
\end{equation*}%
for $x$ $\in $ $G_{m},$ $n$ $\in $ $N.$ Define $I_{n}:=I_{n}\left( 0\right) $
for $n\in N_{+}$. Set $e_{n}:=\left( 0,...,0,1,0,...\right) \in $ $G_{m}$
the $n$ th coordinate of which is $1$ and the rest are zeros $\left( n\in
N\right) .$

If we define the so-called generalized number system based on $m$ in the
following way: $M_{0}:=1,$ $M_{k+1}:=m_{k}M_{k}$ $\left( k\in N\right) ,$
then every $n$ $\in $ $N$ can be uniquely expressed as $n=\sum%
\limits_{j=0}^{\infty }n_{j}M_{j},$ where $n_{j}$ $\in \ Z_{m_{j}}$ $\left(
j\in N_{+}\right) $ and only a finite number of $n_{j}$'s differ from zero.
We also use the following notation: $\left\vert n\right\vert :=$max$\{k\in
N:n_{k}\neq 0\}$ (that is , $M_{|n|}\leq n<M_{|n|+1}$, $n\neq 0$). For every 
$x \in G_{m}$ we denote $\vert x \vert:=\sum\limits_{j=0}^{\infty}\frac{x_j}{%
M_{j+1}}, \left( x_j\text{ }\in Z_{m_{j}}\right)$.

Next, we introduce on $G_{m}$ an orthonormal system, which is called
Vilenkin system. At first define the complex valued functions $r_{k}\left(
x\right) :G_{m}\rightarrow C$, the generalized Rademacher functions, in this
way: 
\begin{equation*}
r_{k}(x):=\exp \frac{2\pi ix_{k}}{m_{k}}\text{\ \ \ }\left( i^{2}=-1,\text{ }%
x\in G_{m},\text{ }k\text{ }\in \text{ }N\right) .
\end{equation*}

Now we define the Vilenkin system $\psi :=\left( \psi _{n}:n\in N\right) $
on $G_{m}$ as follows:

\begin{equation*}
\psi _{n}\left( x\right) :=\prod\limits_{k=0}^{\infty }r_{k}^{n_{k}}\left(
x\right) ,\text{\ \ \ \ \ \ \ }\left( n\text{ }\epsilon \text{ }N\right) .
\end{equation*}

In particular, we call the system the Walsh-Paley system if $m=2.$ Each $%
\psi _{n}$ is a character of $G_{m}$ and all characters of $G_{m}$ are of
this norm. Moreover, $\psi _{n}\left( -x\right) =\bar{\psi}_{n}\left(
x\right) $.

The Dirichlet kernels are defined by

\begin{equation*}
D_{n}:=\sum\limits_{k=0}^{n-1}\psi _{k},\ \ \ \ \ \left( n\text{ }\in
N_{+}\right) .
\end{equation*}

Recall that (see \cite{Gol} or \cite{Sw})

\begin{equation}
\quad \hspace*{0in}D_{M_{n}}\left( x\right) =\left\{ 
\begin{array}{l}
\text{ }M_{n},\text{\thinspace \thinspace \thinspace \thinspace if\thinspace
\thinspace }x\in I_{n}, \\ 
\text{ }0,\text{\thinspace \thinspace \thinspace \thinspace \thinspace if
\thinspace \thinspace }x\notin I_{n}.%
\end{array}%
\right.  \label{for1}
\end{equation}

The Vilenkin system is orthonormal and complete in \ $L^{1}\left(
G_{m}\right) $ ( see \cite{AVDR}).

Next, we introduce some notation with respect to the theory of
two-demonsional Vilenkin system. Let $\tilde{m}$ be a sequence like $m$. The
relation between the sequences $\left( \tilde{m}_{n}\right) $ and \ $\left( 
\tilde{M}_{n}\right) $ is the same as between sequences $\left( m_{n}\right) 
$ and $\left( M_{n}\right) .$ The group $G_{m}\times G_{\tilde{m}}$ is
called a two-dimensional Vilenkin group. The normalized Haar measure is
denoted by $\mu $ as in the one-dimensional case. We also suppose that $m=%
\tilde{m}$ and $G_{m}\times G_{\tilde{m}}=G_{m}^{2}.$

The norm of the space $L^{p}\left( G_{m}^{2}\right) $ is defined by

\begin{equation*}
\left\Vert f\right\Vert _{p}:=\left( \int\limits_{G_{m}^{2}}\left\vert
f\left( x,y\right) \right\vert ^{p}d\mu \left( x,y\right) \right) ^{1/p},\ \
\ \left( 1\leq p<\infty \right) .
\end{equation*}

Denote by $C\left( G_{m}^{2}\right) $ the class of continuous functions on
the group $G_{m}^{2}$, endoved with the supremum norm.

For the sake of brevity in notation, we agree to write $L^{\infty }\left(
G_{m}^{2}\right) $ instead of $C\left( G_{m}^{2}\right) .$

The two-dimensional Fourier coefficients, the rectangular partial sums of
the Fourier series, the Dirichlet kernels with respect to the
two-dimensional Vilenkin system are defined as follows:

\begin{equation*}
\widehat{f}\left( n_{1},n_{2}\right) :=\int\limits_{G_{m}^{2}}f\left(
x,y\right) \bar{\psi}_{n_{1}}\left( x\right) \bar{\psi}_{n_{2}}\left(
y\right) d\mu \left( x,y\right) ,
\end{equation*}

\begin{equation*}
S_{n_{1},n_{2}}\left( x,y,f\right)
:=\sum\limits_{k_{1}=0}^{n_{1}-1}\sum\limits_{k_{2}=0}^{n_{2}-1}\widehat{f}%
\left( k_{1},k_{2}\right) \psi _{k_{1}}\left( x\right) \psi _{k_{2}}\left(
y\right) ,\ 
\end{equation*}

\begin{equation*}
D_{n_{1},n_{2}}\left( x,y\right) :=D_{n_{1}}\left( x\right) D_{n_{2}}\left(
y\right) ,
\end{equation*}%
Denote 
\begin{equation*}
S_{n}^{\left( 1\right) }\left( x,y,f\right) :=\sum\limits_{l=0}^{n-1}%
\widehat{f}\left( l,y\right) \bar\psi _{l}\left( x\right) ,\ 
\end{equation*}

\begin{equation*}
S_{m}^{\left( 2\right) }\left( x,y,f\right) :=\sum\limits_{r=0}^{m-1}%
\widehat{f}\left( x,r\right) \bar\psi _{r}\left( y\right) ,\ 
\end{equation*}
where 
\begin{equation*}
\widehat{f}\left( l,y\right) =\int\limits_{G_{m}}f\left( x,y\right) \psi
_{l}\left( x\right) d\mu \left( x\right)
\end{equation*}
and 
\begin{equation*}
\widehat{f}\left( x,r\right) =\int\limits_{G_{m}}f\left( x,y\right) \psi
_{r}\left( y\right) d\mu \left( y\right) .
\end{equation*}

The $(C,-\alpha )$ means of double Vilenkin-Fourier series are defined as
follows%
\begin{equation*}
\sigma _{n}^{-\alpha }\left( f,x,y\right) =\frac{1}{A_{n-1}^{-\alpha }}%
\sum\limits_{j=1}^{n}A_{n-j}^{-\alpha -1}S_{j,j}\left( f,x,y\right) ,
\end{equation*}

where

\begin{equation*}
A_{0}^{\alpha }=1,\ \ \ \ \ \ \ A_{n}^{\alpha }=\frac{\left( \alpha
+1\right) ...\left( \alpha +n\right) }{n!}.
\end{equation*}

It is well known that (see \cite{Zy})

\begin{equation}
A_{n}^{\alpha }=\sum\limits_{k=0}^{n}A_{k}^{\alpha -1}.\text{ \ \ \ \ \ \ }
\label{for2}
\end{equation}

\begin{equation}
A_{n}^{\alpha }-A_{n-1}^{\alpha }=A_{n}^{\alpha -1}.\ \ \ \ \ \ \ 
\label{for3}
\end{equation}
and 
\begin{equation}
c_{1}(\alpha)n^{\alpha }\leq A_{n}^{\alpha }\leq c_{2} (\alpha)n^{\alpha },\
\ \ \ \ \ \   \label{for4}
\end{equation}
where positive constants $c_{1}$ and $c_{2}$ are dependent on $\alpha$.

The dyadic partial moduli of continuity of a function $f\in L^{p}\left(
G_{m}^{2}\right) \ $\ in the $L^{p}$-norm are defined by

\begin{equation*}
\omega _{1}\left( f,\frac{1}{M_{n}}\right) _{p}=\sup_{u\in I_{n}}\left\Vert
f\left( \cdot +u,\cdot \right) -f\left( \cdot ,\cdot \right) \right\Vert
_{p},
\end{equation*}
and 
\begin{equation*}
\omega _{2}\left( f,\frac{1}{M_{n}}\right) _{p}=\sup_{v\in I_{n}}\left\Vert
f\left( \cdot ,\cdot +v\right) -f\left( \cdot ,\cdot \right) \right\Vert
_{p},
\end{equation*}
while the dyadic mixed modulus of continuity is defined as follows:

\begin{equation*}
\omega _{1,2}\left( f,\frac{1}{M_{n}},\frac{1}{M_{m}}\right) _{p}
\end{equation*}

\begin{equation*}
=\sup_{\left( u,v\right) \in I_{n}\times I_{m}}\left\Vert f\left( \cdot
+u,\cdot +v\right) -f\left( \cdot +u,\cdot \right) -f\left( \cdot ,\cdot
+v\right) +f\left( \cdot ,\cdot \right) \right\Vert _{p}.
\end{equation*}
It is clear that 
\begin{equation*}
\omega _{1,2}\left( f,\frac{1}{M_{n}},\frac{1}{M_{m}}\right) _{p}\leq \omega
_{1}\left( f,\frac{1}{M_{n}}\right) _{p}+\omega _{2}\left( f,\frac{1}{M_{m}}%
\right) _{p}.
\end{equation*}

The dyadic total modulus of continuity is defined by

\begin{equation*}
\omega \left( f,\frac{1}{M_{n}}\right) _{p}=\sup_{\left( u,v\right) \in
I_{n}\times I_{n}}\left\Vert f\left( \cdot +u,\cdot +v\right) -f\left( \cdot
,\cdot \right) \right\Vert _{p}.
\end{equation*}%
The problems of summability of partial sums and Ces\`{a}ro means for
Walsh-Fourier series were studied in \cite{Fi}, \cite{GoAMH}-\cite{Su}, \cite%
{Tev}.

The convergence issue of Fejér (and Ces\`{a}ro ) means on the Walsh and
Vilenkin groups for unbouded case were studies in \cite{GA}-\cite{GATGOG} .

In his monography \cite{Zh} L.V. Zhizhinashvili investigated the behavior of
Ces\`{a}ro $(C, \alpha)-$means for double trigonometric Fourier series in
detail. U.Goginava \cite{GogAn3} studied the analogical question in case of
the Walsh system. In particular, the following theorems were proved:

\begin{Go2}
Let $f$ belong to $L^{p}\left( G_{2}\right) $ for some $p$ $\in $ $\left[
1,\infty \right] $ and \ $\alpha $ $\in $ $\left( 0,1\right) $. Then, for
any \ $2^{k}\leq n<2^{k+1},$ $(k, n\in N)$\ , the inequality

\begin{equation*}
\left\Vert \sigma _{2^{k}}^{-\alpha }\left( f\right) -f\right\Vert _{p}\leq
c\left( \alpha \right) \left\{ 2^{k\alpha }\omega _{1}\left(
f,1/2^{k-1}\right) _{p}+2^{k\alpha }\omega _{2}\left( f,1/2^{k-1}\right)
_{p}+\right.
\end{equation*}

\begin{equation*}
\left. +\sum\limits_{r=0}^{k-2}2^{r-k}\omega _{1}\left( f,1/2^{r}\right)
_{p}+\sum\limits_{s=0}^{k-2}2^{s-k}\omega _{2}\left( f,1/2^{s}\right)
_{p}\right\}
\end{equation*}

holds.
\end{Go2}

\begin{Go3}
\bigskip Let $f$ belong to $L^{p}\left( G_{2}\right) $ for some $p$ $\in $ $%
\left[ 1,\infty \right] $ and \ $\alpha $ $\in $ $\left( 0,1\right) $. Then,
for any \ $2^{k}\leq n<2^{k+1},$ $(k, n\in N)$\ , the inequality

\begin{equation*}
\left\Vert \sigma _{n}^{-\alpha }\left( f\right) -f\right\Vert _{p}\leq
c\left( \alpha \right) \left\{ 2^{k\alpha }k\omega _{1}\left(
f,1/2^{k-1}\right) _{p}+2^{k\alpha }k\omega _{2}\left( f,1/2^{k-1}\right)
_{p}+\right.
\end{equation*}

\begin{equation*}
\left. +\sum\limits_{r=0}^{k-2}2^{r-k}\omega _{1}\left( f,1/2^{r}\right)
_{p}+\sum\limits_{s=0}^{k-2}2^{s-k}\omega _{2}\left( f,1/2^{s}\right)
_{p}\right\}
\end{equation*}

holds.
\end{Go3}

In this paper, we state and prove the analogous results in the case of
double Vilenkin-Fourier series. Our main results read: \newline

\begin{theorem}
\label{T1}Let $f$ belong to $L^{p}\left( G_{m}^{2}\right) $ for some $p$ $%
\in $ $\left[ 1,\infty \right] $ and \ $\alpha $ $\in $ $\left( 0,1\right) $%
. Then, for any \ $M_{k}\leq n<M_{k+1}$ $(k,n\in N)$\ , the inequality
\end{theorem}

\begin{equation*}
\left\Vert \sigma _{M_{k}}^{-\alpha }\left( f\right) -f\right\Vert _{p}\leq
c\left( \alpha \right) \left( \omega _{1}\left( f,1/M_{k-1}\right)
_{p}M_{k}^{\alpha }+\omega _{2}\left( f,1/M_{l-1}\right) _{p}M_{k}^{\alpha
}+\right.
\end{equation*}

\begin{equation*}
\left. +\sum\limits_{r=0}^{k-2}\frac{M_{r}}{M_{k}}\omega _{1}\left(
f,1/M_{r}\right) _{p}+\sum\limits_{s=0}^{k-2}\frac{M_{s}}{M_{k}}\omega
_{2}\left( f,1/M_{s}\right) _{p}\right)
\end{equation*}

holds.

\begin{theorem}
\label{T2}Let $f$ belong to $L^{p}\left( G_{m}^{2}\right) $ for some $p$ $%
\in $ $\left[ 1,\infty \right] $ and \ $\alpha $ $\in $ $\left( 0,1\right) $%
. Then, for any \ $M_{k}\leq n<M_{k+1}$ $(k,n\in N)$\ , the inequality
\end{theorem}

\begin{equation*}
\left\Vert \sigma _{n}^{-\alpha }\left( f\right) -f\right\Vert _{p}\leq
\end{equation*}

\begin{equation*}
c\left( \alpha \right) \left( \omega _{1}\left( f,1/M_{k-1}\right)
_{p}M_{k}^{\alpha }\log n+\omega _{2}\left( f,1/M_{l-1}\right)
_{p}M_{k}^{\alpha }\log n\right.
\end{equation*}

\begin{equation*}
\left. +\sum\limits_{r=0}^{k-2}\frac{M_{r}}{M_{k}}\omega _{1}\left(
f,1/M_{r}\right) _{p}+\sum\limits_{s=0}^{k-2}\frac{M_{s}}{M_{k}}\omega
_{2}\left( f,1/M_{s}\right) _{p}\right)
\end{equation*}

holds. \newline

In order to make the proofs of these Theorems more clear we formulate some
auxiliary Lemmas in Section 2. Some of these Lemmas are new and of
independent interest. The detailed proofs can be found in Section 3. \newpage

\section{AUXILIARY LEMMAS}

In order to prove Theorem 1 and Theorem 2 we need the following Lemmas (see 
\cite{AVDR}, \cite{Gl} and \cite{GogAn4}, respectively)

\begin{lemma}
\label{L0}Let $\alpha _{1}, \alpha _{2},...,\alpha _{n}$ be real numbers.Then
\end{lemma}

\begin{equation*}
\frac{1}{n}\int\limits_{G}\left\vert \sum\limits_{k=1}^{n}\alpha
_{k}D_{k}(x)\right\vert d\mu (x)\leq \frac{c}{\sqrt{n}}\left(
\sum\limits_{k=1}^{n}\alpha _{k}^{2}\right) ^{1/2}.
\end{equation*}

\begin{lemma}
\label{L1} Let $\alpha _{1}, \alpha _{2},...,\alpha _{n}$ be real numbers.
Then

\begin{equation*}
\frac{1}{n}\int\limits_{G_{m}^{2}}\left\vert \sum\limits_{k=1}^{n}\alpha
_{k}D_{k}\left( x\right) D_{k}\left( y\right) \right\vert d\mu \left(
x,y\right) \leq \frac{c}{\sqrt{n}}\left( \sum\limits_{k=1}^{n}\alpha
_{k}^{2}\right) ^{1/2}.
\end{equation*}
\end{lemma}

\begin{lemma}
\label{L2} Let $0\leq j<n_{s}M_{s}$ and $0\leq n_{s}<m_{s}.$ Then

\begin{equation*}
D_{n_{s}M_{s}-j}=D_{n_{s}M_{s}}-\psi _{n_{s}M_{s}-1}\bar{D}_{j}.
\end{equation*}
\end{lemma}

We also need the following new Lemmas of independent interest.

\begin{lemma}
\label{L3} Let $f$ belong to $L^{p}\left( G_{m}^{2}\right) $ for some $p$ $%
\in $ $\left[ 1,\infty \right] $. Then, for every \ $\alpha $ $\in $ $\left(
0,1\right) $, the following inequality holds

\begin{equation*}
I:=\frac{1}{A_{n}^{-\alpha }}\left\Vert
\int\limits_{G_{m}^{2}}\sum\limits_{i=1}^{M_{k-1}}A_{n-i}^{-\alpha
-1}D_{i}\left( u\right) D_{i}\left( v\right) \left[ f\left( \cdot -u,\cdot
-\right) -f\left( \cdot ,\cdot \right) \right] d\mu (u,v)\right\Vert _{p}
\end{equation*}

\begin{equation*}
\left. \leq \ \sum\limits_{r=0}^{k-2}\frac{M_{r}}{M_{k}}\omega _{1}\left(
f,1/M_{r}\right) _{p}+\sum\limits_{s=0}^{k-2}\frac{M_{s}}{M_{k}}\omega
_{2}\left( f,1/M_{s}\right) _{p}\right. ,
\end{equation*}
where $M_{k}\leq n<M_{k+1}.$
\end{lemma}

\begin{lemma}
\label{L4} Let $\alpha $ $\in $ $\left( 0,1\right) $ and $\
p=M_{k},M_{k}+1,....$ Then

\begin{equation*}
II:=\int\limits_{G_{m}^{2}}\left\vert
\sum\limits_{i=1}^{M_{k}}A_{p-i}^{-\alpha -1}D_{i}\left( u\right)
D_{i}\left( v\right) \right\vert d\mu (u,v)\leq c\left( \alpha \right)
<\infty ,\ k=1,2...
\end{equation*}
\end{lemma}

\begin{lemma}
\label{L5} The inequality

\begin{equation*}
III:=\int\limits_{G_{m}^{2}}\left\vert \sum\limits_{i=1}^{n}A_{n-i}^{-\alpha
-1}D_{i}\left( u\right) D_{i}\left( v\right) \right\vert d\mu (u,v)\leq
c\left( \alpha \right) \log n
\end{equation*}

holds.
\end{lemma}

\newpage

\section{The detailed proofs}

\begin{proof}[\textbf{Proof of Lemma 3.}]
Applying Abel's transformation, from (\ref{for2}) we get that

\begin{equation}
I\leq \frac{1}{A_{n}^{-\alpha }}\left\Vert
\int\limits_{G_{m}^{2}}\sum\limits_{i=1}^{M_{k-1}-1}A_{n-i}^{-\alpha
-2}\sum\limits_{l=1}^{i}D_{i}\left( u\right) D_{i}\left( v\right) \left[
f\left( \cdot -u,\cdot -v\right) -f\left( \cdot ,\cdot \right) \right] d\mu
(u,v)\right\Vert _{p}  \label{for5}
\end{equation}

\begin{equation*}
+\frac{1}{A_{n}^{-\alpha }}\left\Vert
\int\limits_{G_{m}^{2}}A_{n-M_{k-1}}^{-\alpha
-1}\sum\limits_{i=1}^{M_{k-1}}D_{i}\left( u\right) D_{i}\left( v\right) %
\left[ f\left( \cdot -u,\cdot -v\right) -f\left( \cdot ,\cdot \right) \right]
d\mu (u,v)\right\Vert _{p}
\end{equation*}

\begin{equation*}
:=I_{1}+I_{2},
\end{equation*}

where the first and the second terms on the right side of inequality (\ref%
{for5}) should be denoted by $I_{1}$ and $I_{2}$ respectively.

For $I_{2}$ we can estimate as follows:

\begin{equation}
I_{2}\leq \frac{1}{A_{n}^{-\alpha }} \biggl\Vert\int%
\limits_{G_{m}^{2}}A_{n-M_{k-1}}^{-\alpha-1}\sum\limits_{r=1}^{k-2}\sum%
\limits_{i=M_{r}}^{M_{r+1}-1}D_{i}\left(u\right) D_{i}\left( v\right)
\label{for6}
\end{equation}

\begin{equation*}
\times \left. \left[ f\left( \cdot -u,\cdot -v\right) -f\left( \cdot ,\cdot
\right) \right]\right.\biggl\Vert_{p}d\mu (u,v)
\end{equation*}

\begin{equation*}
\leq \frac{1}{A_{n}^{-\alpha }}\biggl\Vert\int%
\limits_{G_{m}^{2}}A_{n-M_{k-1}}^{-\alpha-1}\sum\limits_{r=1}^{k-2}\sum%
\limits_{i=M_{r}}^{M_{r+1}-1}D_{i}\left(u\right) D_{i}\left( v\right)
\end{equation*}

\begin{equation*}
\times \left. \left[ f\left( \cdot -u,\cdot -v\right) -S_{M_{r},M_{r}}\left(
\cdot -u,\cdot -v,f\right) \right]\right. d\mu (u,v)\biggl\Vert _{p}
\end{equation*}

\begin{equation*}
+\frac{1}{A_{n}^{-\alpha }}\biggl\Vert \int%
\limits_{G_{m}^{2}}A_{n-M_{k-1}}^{-\alpha
-1}\sum\limits_{r=1}^{k-2}\sum\limits_{i=M_{r}}^{M_{r+1}-1}D_{i}\left(
u\right) D_{i}\left( v\right)
\end{equation*}

\begin{equation*}
\times \left. \left[ S_{M_{r},M_{r}}\left( \cdot -u,\cdot -v,f\right)
-S_{M_{r},M_{r}}\left( \cdot ,\cdot ,f\right) \right] \right.d\mu (u,v)%
\biggl\Vert _{p}
\end{equation*}

\begin{equation*}
+\frac{1}{A_{n}^{-\alpha }}\biggl\Vert \int%
\limits_{G_{m}^{2}}A_{n-M_{k-1}}^{-\alpha
-1}\sum\limits_{r=1}^{k-2}\sum\limits_{i=M_{r}}^{M_{r+1}-1}D_{i}\left(
u\right) D_{i}\left( v\right)
\end{equation*}

\begin{equation*}
\times \left. \left[ S_{M_{r},M_{r}}\left( \cdot ,\cdot ,f\right) -f\left(
\cdot ,\cdot \right) \right] d\mu (u,v)\right. \biggl\Vert%
_{p}:=I_{21}+I_{22}+I_{23},
\end{equation*}%
where the first, the second and the third terms on the right side of
inequality (\ref{for6}) should be denoted by $I_{21}$, $I_{22}$ and $I_{23}$
respectively.

It is evident that

\begin{equation*}
\int\limits_{G_{m}^{2}}\sum\limits_{i=M_{r}}^{M_{r+1}-1}D_{i}\left( u\right)
D_{i}\left( v\right) \left[ S_{M_{r},M_{r}}\left( \cdot -u,\cdot -v,f\right)
-S_{M_{r},M_{r}}\left( \cdot ,\cdot ,f\right) \right] d\mu (u,v)
\end{equation*}

\begin{equation*}
=\sum\limits_{i=M_{r}}^{M_{r+1}-1}\left( \int\limits_{G_{m}^{2}}D_{i}\left(
u\right) D_{i}\left( v\right) S_{M_{r},M_{r}}\left( \cdot -u,\cdot
-v,f\right) d\mu (u,v)-S_{M_{r},M_{r}}\left( \cdot ,\cdot ,f\right) \right)
\end{equation*}

\begin{equation*}
=\sum\limits_{i=M_{r}}^{M_{r+1}-1}\left( S_{i}\left( \cdot ,\cdot
,S_{M_{r},M_{r}}\left( f\right) \right) -S_{M_{r},M_{r}}\left( \cdot ,\cdot
,f\right) \right)
\end{equation*}

\begin{equation*}
=\sum\limits_{i=M_{r}}^{M_{r+1}-1}\left( S_{M_{r},M_{r}}\left( \cdot ,\cdot
,f\right) -S_{M_{r},M_{r}}\left( \cdot ,\cdot ,f\right) \right) =0.
\end{equation*}

Hence,

\begin{equation}
I_{22}=0.  \label{for7}
\end{equation}

Moreover, according to the generalized Minkowski inequality, Lemma \ref{L1}
and by (\ref{for1}) and (\ref{for4}) we obtain that

\begin{equation}
I_{21}\leq \frac{1}{A_{n}^{-\alpha }}\left\vert A_{n-M_{k-1}}^{-\alpha
-1}\right\vert \sum\limits_{r=1}^{k-2}\int\limits_{G_{m}^{2}}\left\vert
\sum\limits_{i=M_{r}}^{M_{r+1}-1}D_{i}\left( u\right) D_{i}\left( v\right)
\right\vert  \label{for8}
\end{equation}

\begin{equation*}
\times \biggl\Vert f\left( \cdot -u,\cdot -v\right) -S_{M_{r},M_{r}}\left(
\cdot -u,\cdot -v,f\right) \biggl\Vert _{p}d\mu (u,v)
\end{equation*}

\begin{equation*}
\leq \frac{c\left( \alpha \right) }{M_{k}}\sum\limits_{r=1}^{k-2}\left(
\omega _{1}\left( f,1/M_{r}\right) _{p}+\omega _{2}\left( f,1/M_{r}\right)
_{p}\right)
\end{equation*}

\begin{equation*}
\times \int\limits_{G_{m}^{2}}\left\vert
\sum\limits_{i=M_{r}}^{M_{r+1}-1}D_{i}\left( x\right) D_{i}\left( y\right)
\right\vert d\mu (u,v)
\end{equation*}

\begin{equation*}
\leq c\left( \alpha \right) \sum\limits_{r=1}^{k-2}\frac{M_{r}}{M_{k}}\left(
\omega _{1}\left( f,1/M_{r}\right) _{p}+\omega _{2}\left( f,1/M_{r}\right)
_{p}\right) .
\end{equation*}
\newline
The estimation of\ $I_{23}$ is analogous to the estimation of $I_{21}$ and
we get that

\begin{equation}
I_{23}\leq c\left( \alpha \right) \sum\limits_{r=1}^{k-2}\frac{M_{r}}{M_{k}}%
\left( \omega _{1}\left( f,1/M_{r}\right) _{p}+\omega _{2}\left(
f,1/M_{r}\right) _{p}\right) .  \label{for9}
\end{equation}%
\newline
Analogously, we can estimate $I_{1}$ in the following way 
\begin{equation}
I_{1}\leq \frac{1}{A_{n}^{-\alpha }}\sum\limits_{r=1}^{k-2}\biggl\Vert %
\int\limits_{G_{m}^{2}}\sum\limits_{i=M_{r}}^{M_{r+1}-1}A_{n-i}^{-\alpha
-2}\sum\limits_{l=1}^{i}D_{l}\left( u\right) D_{l}\left( v\right)
\label{for10}
\end{equation}

\begin{equation*}
\times \left. \left[ f\left( \cdot -u,\cdot -v\right) -S_{M_{r},M_{r}}\left(
\cdot -u,\cdot -v,f\right) \right] d\mu (u,v)\right.\biggl\Vert _{p}
\end{equation*}

\begin{equation*}
+\frac{1}{A_{n}^{-\alpha }}\sum\limits_{r=1}^{k-2}\biggl\Vert %
\int\limits_{G_{m}^{2}}\sum\limits_{i=M_{r}}^{M_{r+1}-1}A_{n-i}^{-\alpha
-2}\sum\limits_{l=1}^{i}D_{l}\left( u\right) D_{l}\left( v\right)
\end{equation*}

\begin{equation*}
\times \left[ S_{M_{r},M_{r}}\left( \cdot -u,\cdot -v,f\right)
-S_{M_{r},M_{r}}\left( \cdot ,\cdot ,f\right) \right] \biggl\Vert _{p}d\mu
(u,v)
\end{equation*}

\begin{equation*}
+\frac{1}{A_{n}^{-\alpha }}\sum\limits_{r=1}^{k-2}\biggl\Vert %
\int\limits_{G_{m}^{2}}\sum\limits_{i=M_{r}}^{M_{r+1}-1}A_{n-i}^{-\alpha
-2}\sum\limits_{l=1}^{i}D_{l}\left( u\right) D_{l}\left( v\right)
\end{equation*}

\begin{equation*}
\times \left. \left[ S_{M_{r},M_{r}}\left( \cdot ,\cdot ,f\right) -f\left(
\cdot ,\cdot \right) \right] d\mu (u,v)\right.\biggl\Vert _{p}
\end{equation*}

\begin{equation*}
\leq \frac{1}{A_{n}^{-\alpha }}\sum\limits_{r=1}^{k-2}\int%
\limits_{G_{m}^{2}}\left\vert
\sum\limits_{i=M_{r}}^{M_{r+1}-1}A_{n-i}^{-\alpha
-2}\sum\limits_{l=1}^{i}D_{l}\left( u\right) D_{l}\left( v\right) \right\vert
\end{equation*}

\begin{equation*}
\times \biggl\Vert f\left( \cdot -u,\cdot -v\right) -S_{M_{r},M_{r}}\left(
\cdot -u,\cdot -v,f\right) \biggl\Vert _{p}d\mu (u,v)
\end{equation*}

\begin{equation*}
+\frac{1}{A_{n}^{-\alpha }}\sum\limits_{r=1}^{k-2}\int\limits_{G_{m}^{2}}%
\left\vert \sum\limits_{i=M_{r}}^{M_{r+1}-1}A_{n-i}^{-\alpha
-2}\sum\limits_{l=1}^{i}D_{l}\left( u\right) D_{l}\left( v\right) \right\vert
\end{equation*}

\begin{equation*}
\times \left. \biggl\Vert S_{M_{r},M_{r}}\left( \cdot ,\cdot ,f\right)
-f\left( \cdot ,\cdot \right)\right. \biggl\Vert _{p} d\mu (u,v)
\end{equation*}

\begin{equation*}
\leq c\left( \alpha \right) M_{k}^{\alpha
}\sum\limits_{r=1}^{k-2}\sum\limits_{i=M_{r}}^{M_{r+1}-1}\left( n-i\right)
^{-\alpha -2}i\left( \omega _{1}\left( f,1/M_{r}\right) _{p}+\omega
_{2}\left( f,1/M_{r}\right) _{p}\right)
\end{equation*}

\begin{equation*}
\leq c\left( \alpha \right) M_{k}^{\alpha
}\sum\limits_{r=1}^{k-2}\sum\limits_{i=M_{r}}^{M_{r+1}-1}\left(
n-M_{r+1}-1\right) ^{-\alpha -2}i\left( \omega _{1}\left( f,1/M_{r}\right)
_{p}+\omega _{2}\left( f,1/M_{r}\right) _{p}\right)
\end{equation*}

\begin{equation*}
\leq c\left( \alpha \right) \sum\limits_{r=0}^{k-2}\frac{M_{r}}{M_{k}}\left(
\omega _{1}\left( f,1/M_{r}\right) _{p}+\omega _{2}\left( f,1/M_{r}\right)
_{p}\right) .
\end{equation*}
\newline
By combining (\ref{for7})-(\ref{for9}) with (\ref{for10}) for $I$ we find
that

\begin{equation}
I\leq c\left( \alpha \right) \sum\limits_{r=0}^{k-2}\frac{M_{r}}{M_{k}}%
\left( \omega _{1}\left( f,1/M_{r}\right) _{p}+\omega _{2}\left(
f,1/M_{r}\right) _{p}\right) .  \label{for12}
\end{equation}
The proof of Lemma 3 is complete.
\end{proof}

\begin{proof}[\textbf{Proof of Lemma 4.}]
It is evident that%
\begin{equation}
II\leq \int\limits_{G_{m}^{2}}\left\vert
\sum\limits_{i=1}^{M_{k}-1}A_{p-M_{k}+i}^{-\alpha -1}D_{M_{k}-i}\left(
u\right) D_{M_{k}-i}\left( v\right) \right\vert d\mu (u,v)  \label{for13}
\end{equation}

\begin{equation*}
+\left\vert A_{p-M_{k}}^{-\alpha -1}\right\vert
\int\limits_{G_{m}^{2}}D_{M_{k}}\left( u\right) D_{M_{k}}\left( v\right)
d\mu (u,v):=II_{1}+II_{2},
\end{equation*}%
where the first and the second terms on the right side of inequality (\ref%
{for13}) should be denoted by $II_{1}$ and $II_{2}$ respectively.

From (\ref{for1}) by $\left\vert A_{p-M_{k}}^{-\alpha -1}\right\vert\leq1$
we get that 
\begin{equation}
II_{2}\leq 1.  \label{for14}
\end{equation}%
Moreover, by Lemma \ref{L2} we have that 
\begin{equation}
II_{1}\leq \int\limits_{G_{m}^{2}}\left\vert
\sum\limits_{i=1}^{M_{k}-1}A_{p-M_{k}+i}^{-\alpha -1}\bar{D}_{i}\left(
u\right) \bar{D}_{i}\left( v\right) \right\vert d\mu (u,v)  \label{for15}
\end{equation}

\begin{equation*}
+\int\limits_{G_{m}^{2}}D_{M_{k}}\left( u\right) \left\vert
\sum\limits_{i=1}^{M_{k}-1}A_{p-M_{k}+i}^{-\alpha -1}\bar{D}_{i}\left(
v\right) \right\vert d\mu (u,v)
\end{equation*}%
\begin{equation*}
+\int\limits_{G_{m}^{2}}D_{M_{k}}\left( v\right) \left\vert
\sum\limits_{i=1}^{M_{k}-1}A_{p-M_{k}+i}^{-\alpha -1}\bar{D}_{i}\left(
u\right) \right\vert d\mu (u,v)
\end{equation*}

\begin{equation*}
+\left\vert \sum\limits_{i=1}^{M_{k}-1}A_{p-M_{k}+i}^{-\alpha -1}\right\vert
\int\limits_{G_{m}^{2}}D_{M_{k}}\left( u\right) D_{M_{k}}\left( v\right)
d\mu (u,v)
\end{equation*}

\begin{equation*}
:=II_{11}+II_{12}+II_{13}+II_{14},
\end{equation*}
where the first, the second, the third and the fourth terms on the right
side of inequality (\ref{for15}) should be denoted by $II_{11}$, $II_{12}$, $%
II_{13}$ and $II_{14}$ respectively.

From (\ref{for1}) and (\ref{for4}) it follows that%
\begin{equation}
II_{14}\leq c\left( \alpha \right) \sum\limits_{v=1}^{\infty }v^{-\alpha
-1}<\infty .  \label{for16}
\end{equation}

By Applying Abel's transformation, in view of Lemma \ref{L1} we have that 
\begin{equation}
II_{11}\leq \int\limits_{G_{m}^{2}}\left\vert
\sum\limits_{i=1}^{M_{k}-2}A_{p-M_{k}+i}^{-\alpha -2}\sum\limits_{l=1}^{i}%
\bar{D}_{l}\left( u\right) \bar{D}_{l}\left( v\right) \right\vert d\mu (u,v)
\label{for17}
\end{equation}

\begin{equation*}
+\int\limits_{G_{m}^{2}}\left\vert A_{p-1}^{-\alpha
-1}\sum\limits_{i=1}^{M_{k}-1}\bar{D}_{i}\left( u\right) \bar{D}_{i}\left(
v\right) \right\vert d\mu (u,v)
\end{equation*}

\begin{equation*}
\leq c\left( \alpha \right) \left\{ \sum\limits_{v=1}^{M_{k}-2}\left(
p-M_{k}+i\right) ^{-\alpha -2}i+\left( p-1\right) ^{-\alpha -1}M_{k}\right\}
\end{equation*}

\begin{equation*}
\leq c\left( \alpha \right) \left\{ \sum\limits_{i=1}^{\infty }i^{-\alpha
-1}+M_{k}^{-\alpha }\right\} <\infty .
\end{equation*}

The estimation of $II_{12}$ and $II_{13}$ are analogous to the estimation of 
$II_{11}$. By Applying Abel's transformation, in view of Lemma \ref{L0} we
find that

\begin{equation}
II_{12}\leq \int\limits_{G_{m}^{2}}D_{M_{k}}\left( u\right) \left\vert
\sum\limits_{i=1}^{M_{k}-2}A_{p-M_{k}+i}^{-\alpha -2}\sum\limits_{l=1}^{i}%
\bar{D}_{l}\left( v\right) \right\vert d\mu (u,v)  \label{for18}
\end{equation}

\begin{equation*}
+\int\limits_{G_{m}^{2}}D_{M_{k}}\left( u\right) \left\vert A_{p-1}^{-\alpha
-1}\sum\limits_{i=1}^{M_{k}-1}\bar{D}_{i}\left( v\right) \right\vert d\mu
(u,v)
\end{equation*}

\begin{equation*}
\leq c\left( \alpha \right) \left\{ \sum\limits_{v=1}^{M_{k}-2}\left(
p-M_{k}+i\right) ^{-\alpha -2}i+\left( p-1\right) ^{-\alpha -1}M_{k}\right\}
\end{equation*}

\begin{equation*}
\leq c\left( \alpha \right) \left\{ \sum\limits_{i=1}^{\infty }i^{-\alpha
-1}+M_{k}^{-\alpha }\right\} <\infty .
\end{equation*}

and

\begin{equation}
III_{12}\leq \int\limits_{G_{m}^{2}}D_{M_{k}}\left( v\right) \left\vert
\sum\limits_{i=1}^{M_{k}-2}A_{p-M_{k}+i}^{-\alpha -2}\sum\limits_{l=1}^{i}%
\bar{D}_{l}\left( u\right) \right\vert d\mu (u,v)  \label{for19}
\end{equation}

\begin{equation*}
+\int\limits_{G_{m}^{2}}D_{M_{k}}\left( v\right) \left\vert A_{p-1}^{-\alpha
-1}\sum\limits_{i=1}^{M_{k}-1}\bar{D}_{i}\left( u\right) \right\vert d\mu
(u,v)
\end{equation*}

\begin{equation*}
\leq c\left( \alpha \right) \left\{ \sum\limits_{v=1}^{M_{k}-2}\left(
p-M_{k}+i\right) ^{-\alpha -2}i+\left( p-1\right) ^{-\alpha -1}M_{k}\right\}
\end{equation*}

\begin{equation*}
\leq c\left( \alpha \right) \left\{ \sum\limits_{i=1}^{\infty }i^{-\alpha
-1}+M_{k}^{-\alpha }\right\} <\infty .
\end{equation*}
\newline
The proof is complete by combining (\ref{for13})-(\ref{for19}).
\end{proof}

\begin{proof}[\textbf{Proof of Lemma 5.}]
Let 
\begin{equation*}
n=n_{k_{1}}M_{k_{1}}+...+n_{k_{s}}M_{k_{s}},\ k_{1}>...>k_{s}\geq 0.
\end{equation*}

Denote

\begin{equation*}
n^{\left( i\right) }=n_{k_{i}}M_{k_{i}}+...+n_{k_{s}}M_{k_{s}},\text{ \ \ }%
i=1,2,...s.
\end{equation*}

Since ( see \cite{Gol}) 
\begin{equation}
D_{j+n_{A}M_{A}}=D_{n_{A}M_{A}}+\psi _{n_{A}M_{A}}D_{j},  \label{for20}
\end{equation}
we find that

\begin{equation}
III\leq \int\limits_{G_{m}^{2}}\left\vert
\sum\limits_{i=1}^{n_{k_{1}}M_{k_{1}}}A_{n-i}^{-\alpha -1}D_{i}\left(
u\right) D_{i}\left( v\right) \right\vert d\mu (u,v)  \label{for21}
\end{equation}

\begin{equation*}
+\int\limits_{G_{m}^{2}}\left\vert \sum\limits_{i=1}^{n^{\left( 2\right)
}}A_{n^{\left( 2\right) }-i}^{-\alpha -1}D_{i}\left( u\right) D_{i}\left(
v\right) \right\vert d\mu (u,v)
\end{equation*}

\begin{equation*}
+\int\limits_{G_{m}^{2}}D_{n_{k_{1}}M_{k_{1}}}\left( u\right)
D_{n_{k_{1}}M_{k_{1}}}\left( v\right) \left\vert
\sum\limits_{i=1}^{n^{\left( 2\right) }}A_{n^{\left( 2\right) }-i}^{-\alpha
-1}\right\vert d\mu (u,v)
\end{equation*}

\begin{equation*}
+\int\limits_{G_{m}^{2}}D_{n_{k_{1}}M_{k_{1}}}\left( u\right) \left\vert
\sum\limits_{i=1}^{n^{\left( 2\right) }}A_{n^{\left( 2\right) }-i}^{-\alpha
-1}D_{i}\left( v\right) \right\vert d\mu (u,v)
\end{equation*}

\begin{equation*}
+\int\limits_{G_{m}^{2}}D_{n_{k_{1}}M_{k_{1}}}\left( v\right) \left\vert
\sum\limits_{i=1}^{n^{\left( 2\right) }}A_{n^{\left( 2\right) }-i}^{-\alpha
-1}D_{i}\left( u\right) \right\vert d\mu (u,v)
\end{equation*}

\begin{equation*}
:=III_{1}+III_{2}+III_{3}+III_{4}+III_{5},
\end{equation*}
where the first, the second, the third, the fourth and the fifth terms on
the right side of inequality (\ref{for21}) should be denoted by $III_{1}$, $%
III_{2}$, $III_{3}$, $III_{4}$ and $III_{5}$ respectively.

By (\ref{for1}) we have that%
\begin{equation}
III_{3}\leq c\left( \alpha \right) .  \label{for22}
\end{equation}

Moreover, since (see \cite{Sh}  )

\begin{equation}
\left\vert \sum\limits_{i=1}^{n}A_{n-i}^{-\alpha -1}D_{i}\left( u\right)
\right\vert =O\left(\vert u \vert^{\alpha -1}\right) .  \label{for23}
\end{equation}

for $III_{4}$ we get that

\begin{equation}
III_{4}\leq \int\limits_{G_{m}^{2}}D_{n_{k_{1}}M_{k_{1}}}\left( u\right)
|v|^{\alpha -1}d\mu (u,v)  \label{for24}
\end{equation}

\begin{equation*}
\leq \int\limits_{G_{m}}|v|^{\alpha -1}d\mu \left( v\right) =\frac{1}{\alpha 
}<\infty .
\end{equation*}

Analogously, we find that

\begin{equation}
III_{5}\leq \int\limits_{G_{m}^{2}}D_{n_{k_{1}}M_{k_{1}}}\left( v\right)
|u|^{\alpha -1}d\mu (u,v)  \label{for25}
\end{equation}

\begin{equation*}
\leq \int\limits_{G_{m}}|u|^{\alpha -1}d\mu \left( v\right) =\frac{1}{\alpha 
}<\infty .
\end{equation*}

For $r\in \{0,...m_{A}-1\},$ $0\leq \ j<M_{A}$ , ( see \cite{Gol}), it
yields that

\begin{equation*}
D_{j+rM_{A}}=\left( \sum\limits_{q=0}^{r-1}\psi _{M_{A}}^{q}\right)
D_{M_{A}}+\psi _{M_{A}}^{r}D_{j}.
\end{equation*}

Thus, we have that

\begin{equation*}
\int\limits_{G_{m}^{2}}\sum\limits_{i=1}^{n_{k_{1}}M_{k_{1}}-1}A_{n-i}^{-%
\alpha -1}D_{i}\left( u\right) D_{i}\left( v\right) d\mu (u,v)
\end{equation*}

\begin{equation*}
\leq
\int\limits_{G_{m}^{2}}\sum\limits_{r=0}^{n_{k_{1}}-1}\sum%
\limits_{i=0}^{M_{k_{1}}-1}A_{n-i-rM_{k_{1}}}^{-\alpha
-1}D_{i+rM_{k_{1}}}\left( u\right) D_{i+rM_{k_{1}}}\left( v\right) d\mu (u,v)
\end{equation*}

\begin{equation*}
\leq
\int\limits_{G_{m}^{2}}\sum\limits_{r=0}^{n_{k_{1}}-1}\sum%
\limits_{i=0}^{M_{k_{1}}-1}A_{n-i-rM_{k_{1}}}^{-\alpha -1}\left(
\sum\limits_{q=0}^{r-1}\psi _{M_{k_{1}}}^{q}\right) D_{M_{k_{1}}}\left(
u\right)
\end{equation*}

\begin{equation*}
\times \left( \sum\limits_{q=0}^{r-1}\psi _{M_{k_{1}}}^{q}\right)
D_{M_{k_{1}}}\left( v\right) d\mu (u,v)
\end{equation*}

\begin{equation*}
+\int\limits_{G_{m}^{2}}\sum\limits_{r=0}^{n_{k_{1}}-1}\sum%
\limits_{i=0}^{M_{k_{1}}-1}A_{n-i-rM_{k_{1}}}^{-\alpha -1}\left(
\sum\limits_{q=0}^{r-1}\psi _{M_{k_{1}}}^{q}\right) D_{M_{k_{1}}}\left(
u\right) \psi _{M_{A}}^{r}D_{i}\left( v\right) d\mu (u,v)
\end{equation*}

\begin{equation*}
+\int\limits_{G_{m}^{2}}\sum\limits_{r=0}^{n_{k_{1}}-1}\sum%
\limits_{i=0}^{M_{k_{1}}-1}A_{n-i-rM_{k_{1}}}^{-\alpha -1}\psi
_{M_{A}}^{r}D_{i}\left( u\right) \left( \sum\limits_{q=0}^{r-1}\psi
_{M_{k_{1}}}^{q}\right) D_{M_{k_{1}}}\left( v\right) d\mu (u,v)
\end{equation*}
\newline
\begin{equation*}
+\int\limits_{G_{m}^{2}}\sum\limits_{r=0}^{n_{k_{1}}-1}\sum%
\limits_{i=0}^{M_{k_{1}}-1}A_{n-i-rM_{k_{1}}}^{-\alpha -1}\psi
_{M_{A}}^{r}D_{i}\left( u\right) \psi _{M_{A}}^{r}D_{i}\left( v\right) d\mu
(u,v).
\end{equation*}

On the other hand, by (\ref{for1}) and (\ref{for4}) we obtain that

\begin{equation*}
\int\limits_{G_{m}^{2}}A_{n-n_{k_{1}}M_{k_{1}}}^{-\alpha
-1}D_{n_{k_{1}}M_{k_{1}}}\left( u\right) D_{n_{k_{1}}M_{k_{1}}}\left(
v\right) d\mu (u,v)\leq c\left( \alpha \right).
\end{equation*}
\newline

Consequently, for $III_{1}$ we have the estimate \newline
\begin{equation}
III_{1}\leq \int\limits_{G_{m}^{2}}D_{M_{k_{1}}}\left( u\right)
D_{M_{k_{1}}}\left( v\right) \left\vert
\sum\limits_{r=0}^{n_{k_{1}}-1}\sum%
\limits_{i=1}^{M_{k_{1}}}A_{n-i-rM_{k_{1}}}^{-\alpha -1}\right\vert d\mu
(u,v)  \label{for25.1}
\end{equation}

\begin{equation*}
+\int\limits_{G_{m}^{2}}D_{M_{k_{1}}}\left( u\right) \left\vert
\sum\limits_{r=0}^{n_{k_{1}}-1}\sum%
\limits_{i=1}^{M_{k_{1}}}A_{n-i-rM_{k_{1}}}^{-\alpha -1}D_{i}\left( v\right)
\right\vert d\mu (u,v)
\end{equation*}

\begin{equation*}
+\int\limits_{G_{m}^{2}}D_{M_{k_{1}}}\left( v\right) \left\vert
\sum\limits_{r=0}^{n_{k_{1}}-1}\sum%
\limits_{i=1}^{M_{k_{1}}}A_{n-i-rM_{k_{1}}}^{-\alpha -1}D_{i}\left( u\right)
\right\vert d\mu (u,v)
\end{equation*}

\begin{equation*}
+\int\limits_{G_{m}^{2}}\left\vert
\sum\limits_{r=0}^{n_{k_{1}}-1}\sum%
\limits_{i=1}^{M_{k_{1}}}A_{n-i-rM_{k_{1}}}^{-\alpha -1}D_{i}\left( u\right)
D_{i}\left( v\right) \right\vert d\mu (u,v)+c\left( \alpha \right)
\end{equation*}

\begin{equation*}
:=III_{11}+III_{12}+III_{13}+III_{14}+c\left( \alpha \right) .
\end{equation*}
where the first, the second, the third and the fourth terms on the right
side of inequality (\ref{for25.1}) should be denoted by $III_{11}$, $%
III_{12} $, $III_{13}$ and $III_{14}$ respectively.

From Lemma 4 we have that

\begin{equation}
III_{14}\leq c\left( \alpha \right) .  \label{for25.2}
\end{equation}

The estimation of\ $III_{11}$ is analogous to the estimation of $III_{3}$
and we find that 
\begin{equation}
III_{11}\leq c\left( \alpha \right) .  \label{for25.3}
\end{equation}

The estimation of\ $III_{12}$ and $III_{13}$ is analogous to the estimation
of $III_{4}$ and we obtain that 
\begin{equation}
III_{12}< \infty ,  \label{for25.4}
\end{equation}
and 
\begin{equation}
III_{13}< \infty .  \label{for25.5}
\end{equation}

After substituting (\ref{for22}) and (\ref{for24})- (\ref{for25.5}) into (%
\ref{for21}) we conclude that

\begin{equation*}
\int\limits_{G_{m}^{2}}\left\vert \sum\limits_{i=1}^{n}A_{n-i}^{-\alpha
-1}D_{i}\left( u\right) D_{i}\left( v\right) \right\vert d\mu (u,v)
\end{equation*}

\begin{equation*}
\leq \int\limits_{G_{m}^{2}}\left\vert \sum\limits_{i=1}^{n^{\left( 2\right)
}}A_{n^{\left( 2\right) }-i}^{-\alpha -1}D_{i}\left( u\right) D_{i}\left(
v\right) \right\vert d\mu (u,v)+c\left( \alpha \right)
\end{equation*}

\begin{equation*}
\leq ...\leq \int\limits_{G_{m}^{2}}\left\vert \sum\limits_{i=1}^{n^{\left(
s\right) }}A_{n^{\left( s\right) }-i}^{-\alpha -1}D_{i}\left( u\right)
D_{i}\left( v\right) \right\vert d\mu (u,v)+c\left( \alpha \right) s
\end{equation*}

\begin{equation*}
\leq c\left( \alpha \right) +c\left( \alpha \right) s\leq c\left( \alpha
\right) \log n.
\end{equation*}

The proof is complete.
\end{proof}

Now we are ready to prove the main results

\begin{proof}[\textbf{Proof of Theorem 1.}]
It is evident that%
\begin{equation}
\left\Vert \sigma _{M_{k}}^{-\alpha }\left( f\right) -f\right\Vert _{p}
\label{for27}
\end{equation}

\begin{equation*}
\leq \frac{1}{A_{M_{k}-1}^{-\alpha }}\left\Vert
\int\limits_{G_{m}^{2}}\sum\limits_{i=1}^{M_{k-1}}A_{M_{k}-i}^{-\alpha
-1}D_{i}\left( u\right) D_{i}\left( v\right) \left[ f\left( \cdot -u,\cdot
-v\right) -f\left( \cdot ,\cdot \right) \right] d\mu (u,v)\right\Vert _{p}
\end{equation*}

\begin{equation*}
+\frac{1}{A_{M_{k}-1}^{-\alpha }}\left\Vert
\int\limits_{G_{m}^{2}}\sum\limits_{i=M_{k-1}+1}^{M_{k}}A_{M_{k}-i}^{-\alpha
-1}D_{i}\left( u\right) D_{i}\left( v\right) \left[ f\left( \cdot -u,\cdot
-v\right) -f\left( \cdot ,\cdot \right) \right] d\mu (u,v)\right\Vert _{p}
\end{equation*}

\begin{equation*}
:=I+II.
\end{equation*}

From Lemma \ref{L4} it follows that

\begin{equation}
I\leq c\left( \alpha \right) \sum\limits_{r=0}^{k-2}\frac{M_{r}}{M_{k}}%
\left( \omega _{1}\left( f,1/M_{r}\right) _{p}+\omega _{2}\left(
f,1/M_{r}\right) _{p}\right) .  \label{for28}
\end{equation}

Moreover, for II we have the estimate

\begin{equation}
II\leq \frac{1}{A_{M_{k}-1}^{-\alpha }}\biggl\Vert\int\limits_{G_{m}^{2}}%
\sum\limits_{i=M_{k-1}+1}^{M_{k}}A_{M_{k}-i}^{-\alpha-1}D_{i}\left( u\right)
D_{i}\left( v\right)  \label{for29}
\end{equation}

\begin{equation*}
\times \left. \left[ f\left( \cdot -u,\cdot -v\right) -S_{M_{k-1}}^{\left(
1\right) }\left( \cdot -u,\cdot -v,f\right) \right] d\mu (u,v)\right.%
\biggl\Vert _{p}
\end{equation*}

\begin{equation*}
+\frac{1}{A_{M_{k}-1}^{-\alpha }}\biggl\Vert \int\limits_{G_{m}^{2}}\sum%
\limits_{i=M_{k-1}+1}^{M_{k}}A_{M_{k}-i}^{-\alpha -1}D_{i}\left( u\right)
D_{i}\left( v\right)
\end{equation*}

\begin{equation*}
\times \left. \left[ S_{M_{k-1}}^{\left( 1\right) }\left( \cdot -u,\cdot
-v,f\right) -f\left( \cdot ,\cdot \right) \right] d\mu (u,v) \right.%
\biggl\Vert _{p}:=II_{1}+II_{2},
\end{equation*}
where the first and the second terms on the right side of inequality (\ref%
{for29}) should be denoted by $II_{1}$ and $II_{2}$ respectively.

In view of generalized Minkowski inequality, by (\ref{for4}) and using Lemma %
\ref{L4} we get that

\begin{equation}
II_{1}\leq \frac{1}{A_{M_{k}-1}^{-\alpha }}\int\limits_{G_{m}^{2}}\left\vert
\sum\limits_{i=M_{k-1}+1}^{M_{k}}A_{M_{k}-i}^{-\alpha -1}D_{i}\left(
u\right) D_{i}\left( v\right) \right\vert  \label{for30}
\end{equation}

\begin{equation*}
\times \left\Vert f\left( \cdot -u,\cdot -v\right) -S_{M_{k-1}}^{\left(
1\right) }\left( \cdot -u,\cdot -v,f\right) \right\Vert _{p}d\mu (u,v)
\end{equation*}

\begin{equation*}
\leq c\left( \alpha \right) M_{k}^{\alpha }\omega _{1}\left(
f,1/M_{k-1}\right) _{p}.
\end{equation*}

The estimation of\ $II_{2}$ is analogous to the estimation of $II_{1}$ and
we find that 
\begin{equation}
II_{2}\leq c\left( \alpha \right) M_{k}^{\alpha }\omega _{2}\left(
f,1/M_{k-1}\right) _{p}.  \label{for31}
\end{equation}

Combining (\ref{for27})- (\ref{for31}) we receive the proof of Theorem 1.
\end{proof}

\begin{proof}[\textbf{Proof of Theorem 2.}]
It is evident that%
\begin{equation}
\left\Vert \sigma _{n}^{-\alpha }\left( f\right) -f\right\Vert _{p}
\label{for32}
\end{equation}%
\begin{equation*}
\leq \frac{1}{A_{n-1}^{-\alpha }}\left\Vert
\int\limits_{G_{m}^{2}}\sum\limits_{i=1}^{M_{k-1}}A_{n-i}^{-\alpha
-1}D_{i}\left( u\right) D_{i}\left( v\right) \left[ f\left( \cdot -u,\cdot
-v\right) -f\left( \cdot ,\cdot \right) \right] d\mu (u,v)\right\Vert _{p}
\end{equation*}

\begin{equation*}
+\frac{1}{A_{n-1}^{-\alpha }}\left\Vert
\int\limits_{G_{m}^{2}}\sum\limits_{i=M_{k-1}+1}^{M_{k}}A_{n-i}^{-\alpha
-1}D_{i}\left( u\right) D_{i}\left( v\right) \left[ f\left( \cdot -u,\cdot
-v\right) -f\left( \cdot ,\cdot \right) \right] d\mu (u,v)\right\Vert _{p}
\end{equation*}

\begin{equation*}
+\frac{1}{A_{n-1}^{-\alpha }}\left\Vert
\int\limits_{G_{m}^{2}}\sum\limits_{i=M_{k}+1}^{n}A_{n-i}^{-\alpha
-1}D_{i}\left( u\right) D_{i}\left( v\right) \left[ f\left( \cdot -u,\cdot
-v\right) -f\left( \cdot ,\cdot \right) \right] d\mu (u,v)\right\Vert _{p}
\end{equation*}

\begin{equation*}
:=I+II+III,
\end{equation*}
where the first, the second and the third terms on the right side of
inequality (\ref{for32}) should be denoted by $I$, $II$ and $III$
respectively.

From Lemma \ref{L3} it follows that 
\begin{equation}
I\leq \ c\left( \alpha \right) \sum\limits_{r=0}^{k-2}\frac{M_{r}}{M_{k}}%
\left( \omega _{1}\left( f,1/M_{r}\right) _{p}+\omega _{2}\left(
f,1/M_{r}\right) _{p}\right) .  \label{for33}
\end{equation}

Next, we repeat the arguments just in the same way as in the proof of
Theorem 1 and find that \newline
\begin{equation}
II\leq \ c\left( \alpha \right) M_{k}^{\alpha }\left( \omega _{1}\left(
f,1/M_{k-1}\right) _{p}+\omega _{2}\left( f,1/M_{k-1}\right) _{p}\right) .
\label{for34}
\end{equation}

On the other hand, for III we have that

\begin{equation}
III\leq \frac{1}{A_{n-1}^{-\alpha }}\biggl\Vert \int\limits_{G_{m}^{2}}\sum%
\limits_{i=M_{k}+1}^{n}A_{n-i}^{-\alpha -1}D_{i}\left( u\right) D_{i}\left(
v\right)  \label{for35}
\end{equation}

\begin{equation*}
\times \left. \left[ f\left( \cdot -u,\cdot -v\right) -f\left( \cdot ,\cdot
\right) \right] \right. \biggl\Vert _{p}d\mu (u,v)
\end{equation*}

\begin{equation*}
\leq \frac{1}{A_{n}^{-\alpha }}\biggl\Vert \int\limits_{G_{m}^{2}}\sum%
\limits_{i=M_{k}+1}^{n}A_{n-i}^{-\alpha -1}D_{i}\left( u\right) D_{i}\left(
v\right) \times
\end{equation*}

\begin{equation*}
\times \left. \left[ f\left( \cdot -u,\cdot -v\right) -S_{M_{k},M_{k}}\left(
\cdot -u,\cdot -v,f\right) \right] d\mu (u,v) \right.\biggl\Vert _{p}
\end{equation*}

\begin{equation*}
\leq \frac{1}{A_{n}^{-\alpha }}\biggl\Vert \int\limits_{G_{m}^{2}}\sum%
\limits_{i=M_{k}+1}^{n}A_{n-i}^{-\alpha -1}D_{i}\left( u\right) D_{i}\left(
v\right) \times
\end{equation*}

\begin{equation*}
\times \left. \left[ S_{M_{k},M_{k}}\left( \cdot -u,\cdot -v,f\right)
-S_{M_{k},M_{k}}\left( \cdot ,\cdot ,f\right) \right] \right. d\mu (u,v)%
\biggl\Vert _{p}
\end{equation*}

\begin{equation*}
\leq \frac{1}{A_{n}^{-\alpha }}\biggl\Vert \int\limits_{G_{m}^{2}}\sum%
\limits_{i=M_{k}+1}^{n}A_{n-i}^{-\alpha -1}D_{i}\left( u\right) D_{i}\left(
v\right) \times
\end{equation*}

\begin{equation*}
\times \left. \left[ S_{M_{k},M_{k}}\left( \cdot ,\cdot ,f\right) -f\left(
\cdot ,\cdot \right) \right] \right. d\mu (u,v)\biggl\Vert %
_{p}:=III_{1}+III_{2}+III_{3},
\end{equation*}%
where the first, the second and the third terms on the right side of
inequality (\ref{for35}) should be denoted by $III_{1}$, $III_{2}$ and $%
III_{3}$ respectively.

It is easy to show that%
\begin{equation}
III_{2}=0.  \label{for36}
\end{equation}%
\newline
According to the generalized Minkowski inequality and by using Lemma 5 for $%
III_{1}$ we obtain that%
\begin{equation}
III_{1}\leq \frac{1}{A_{n}^{-\alpha }}\int\limits_{G_{m}^{2}}\left\vert
\sum\limits_{i=M_{k}+1}^{n}A_{n-i}^{-\alpha -1}D_{i}\left( u\right)
D_{i}\left( v\right) \right\vert  \label{for37}
\end{equation}

\begin{equation*}
\times \biggl\Vert f\left( \cdot -u,\cdot -v\right) -S_{M_{r},M_{r}}\left(
\cdot -u,\cdot -v,f\right) \biggl\Vert _{p}d\mu (u,v)
\end{equation*}
\newline
\begin{equation*}
\leq c\left( \alpha \right) M_{k}^{\alpha }\left( \omega _{1}\left(
f,1/M_{k-1}\right) _{p}+\omega _{2}\left( f,1/M_{k-1}\right) _{p}\right)
\end{equation*}
\newline
\begin{equation*}
\times \int\limits_{G_{m}^{2}}\left\vert
\sum\limits_{v=M_{k}+1}^{n}A_{n-v}^{-\alpha -1}D_{v}\left( u\right)
D_{v}\left( v\right) \right\vert d\mu (u,v)\newline
\end{equation*}
\newline
\begin{equation*}
\leq c\left( \alpha \right) M_{k}^{\alpha }\log n\left( \omega _{1}\left(
f,1/M_{k-1}\right) _{p}+\omega _{2}\left( f,1/M_{k-1}\right) _{p}\right) .
\end{equation*}
\newline
The estimation of $III_{3}$ is analogous to the estimation of $III_{2}$ and
we find that

\begin{equation}
III_{3}\leq c\left( \alpha \right) M_{k}^{\alpha }\log n\left( \omega
_{1}\left( f,1/M_{k-1}\right) _{p}+\omega _{2}\left( f,1/M_{k-1}\right)
_{p}\right) .  \label{for38}
\end{equation}

After substituting (\ref{for33})- (\ref{for34}), (\ref{for38}) into (\ref%
{for32}), we receive the proof of Theorem 2.
\end{proof}

\textbf{Author details}

$^{1}$ The Artic University of Norway, Campus Narvik, P.O. Box 385, N-8505, Narvik, Norway.

$^{2}$ The Artic University of Norway, Campus Narvik, P.O. Box 385, N-8505, Narvik, Norway.

\textbf{Authors' contributions}

The authors contributed equally to the writing of this paper. Both authors
approved the final version of the manuscript.

\textbf{Competing interests}

The authors declare that they have no competing interests.

\textbf{Acknowledgements}

The authors would like to thank the referees for helpful suggestions.

\end{document}